\documentclass[amssymb,aps,floatfix,showkeys,showpacs,twocolumn]{revtex4-1}

\usepackage{graphicx}
\usepackage{amsmath}
\usepackage[pdfauthor="Richard J. Mathar",colorlinks=true,bookmarks=false,citecolor=blue]{hyperref}

\begin{document}

\title[RiemCirc: Nodes and Weights of Integration on the Circle]{RiemCirc: A Generator of Nodes and Weights for Riemann Integration on the Circle}

\author{Richard J. Mathar}
\homepage{https://www.mpia-hd.mpg.de/~mathar}
\email{mathar@mpia-hd.mpg.de}
\affiliation{Max-Planck Institute of Astronomy, K\"onigstuhl 17, 69117 Heidelberg, Germany}


\date{\today}

\begin{abstract}
\texttt{RiemCirc} is a C++ program which allocates points inside the unit circle
for numerical quadrature on the circle, aiming at homogeneous
equidistant distribution. The weights of the quadrature rule are computed
by the area of the tiles that surround these nodes. The shapes of the
areas are polygonal, defined
by Voronoi tessellation.
\end{abstract}

\maketitle
\section{Scope} 
\subsection{Quadrature}

Numerical integration of functions defined on the unit circle
is generally discussed with a wish to preserve the (or some form
of subgroup) symmetry
of the circle while choosing a set of nodes and weights
of sampling the area
\cite{HeoJCAM112,CoolsJCAM20,VerlindenNM61}.
The integral is canonically replaced by a finite Riemann sum
over $N$ nodes with radial and azimuthal coordinates $r_j$ and $\varphi_j$,
\begin{equation}
\int f(r,\varphi) d^2r \approx \sum_{j=1}^N w_j f(r_j,\varphi_j)
.
\end{equation}
The standard procedure is to split the integral into a product
of integrals over $r$ and $\varphi$ and to select
separate one-dimensional Gaussian integration formulas
and associated nodal weights
\cite{StroudMathComp17,TrefethenSIAMR50,BojanovNM80}.
(We will not be interested in integrations along the
perimeter of the circle which typically arises in integrations
in the complex plane
\cite{DaruisJCAM140,HunterBIT35}.)

The disadvantage of the product representation is
a clumpy structure of the nodes in the circle if
$N$ is not large. The program that is presented here
ignores aspects of symmetry  conservation and attempts to
distribute $N$ points well-balanced in the sense that
the weights $w_j$ have approximately equal values, in
fact, minimizing their variance.

\subsection{Voronoi Tessellation}
The idea is that the distances between the nodes
are equilibrated by starting from a randomized set
of locations, then moving the nodes iteratively as if 
some repelling force was adjusting their places.
This is put into concrete by (i) subdividing the area
inside the circle into tiles around the nodes
akin to a Voronoi tessellation, and (ii) moving the individual
nodes to the center-of-mass of their tile, and iterating
these two steps a few times to achieve some stable
configuration.

\begin{figure}
\includegraphics[scale=0.4]{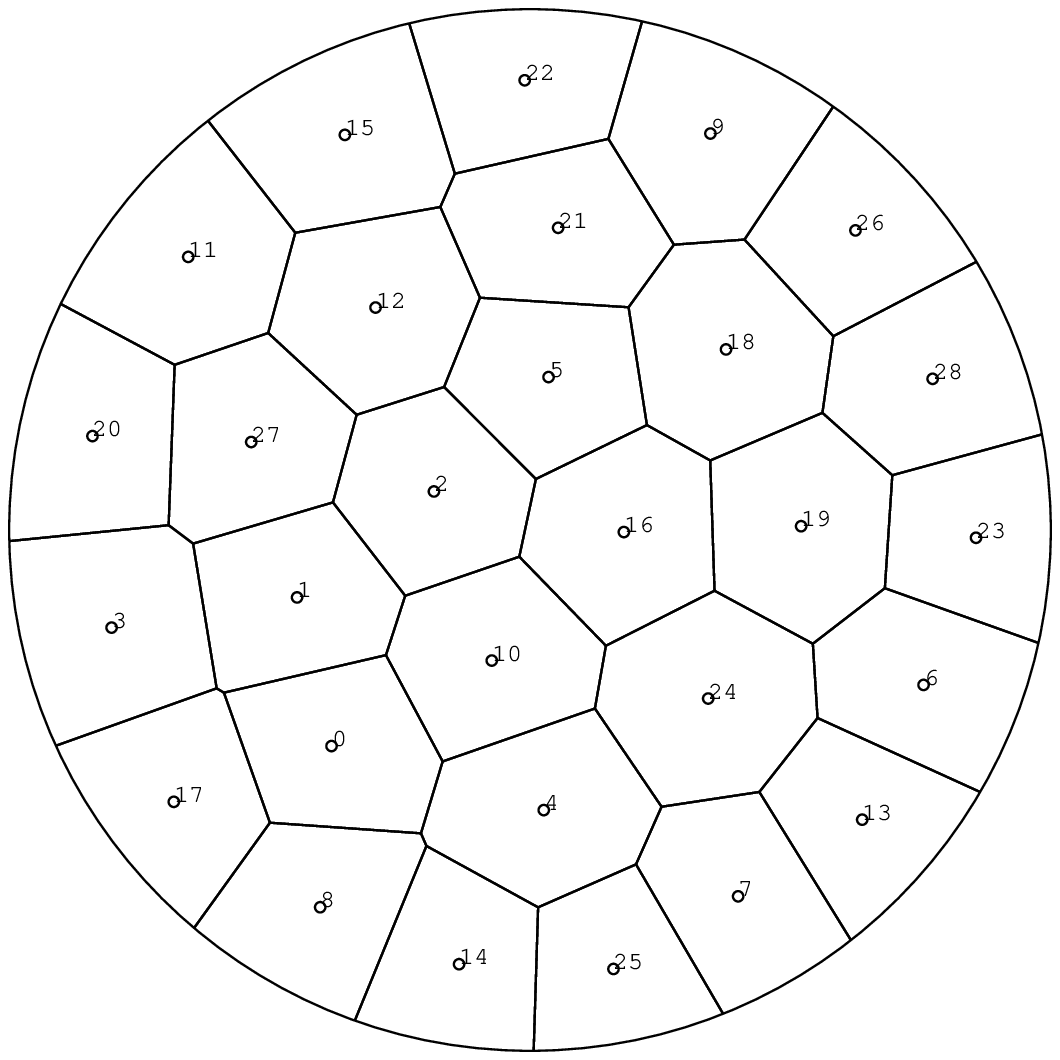}
\caption{
Illustration of a result on $29$ nodes.
}
\label{fig.ex1}
\end{figure}

\begin{figure}
\includegraphics[scale=0.4]{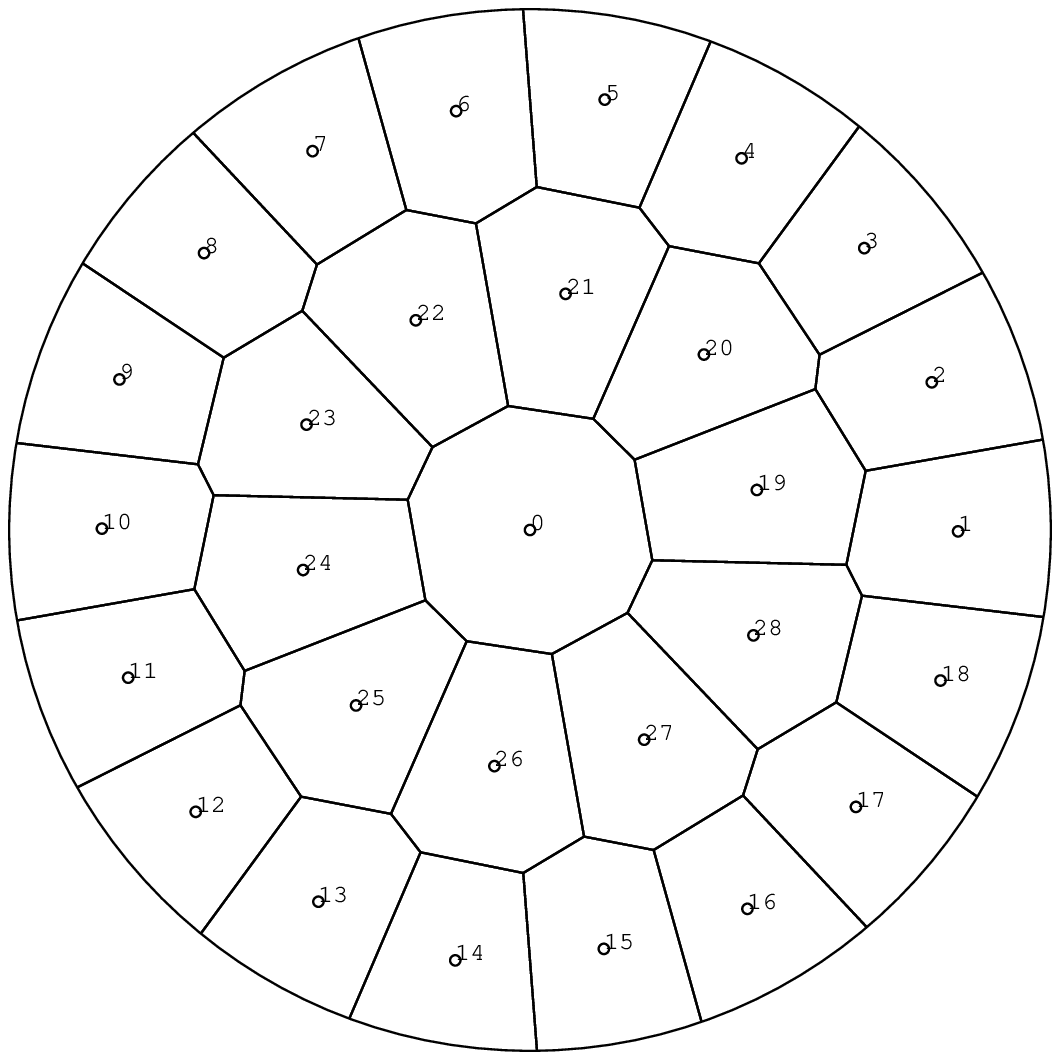}
\caption{
Illustration of a ``centered'' result on $29$ nodes.
}
\label{fig.ex2}
\end{figure}

\begin{figure}
\includegraphics[scale=0.4]{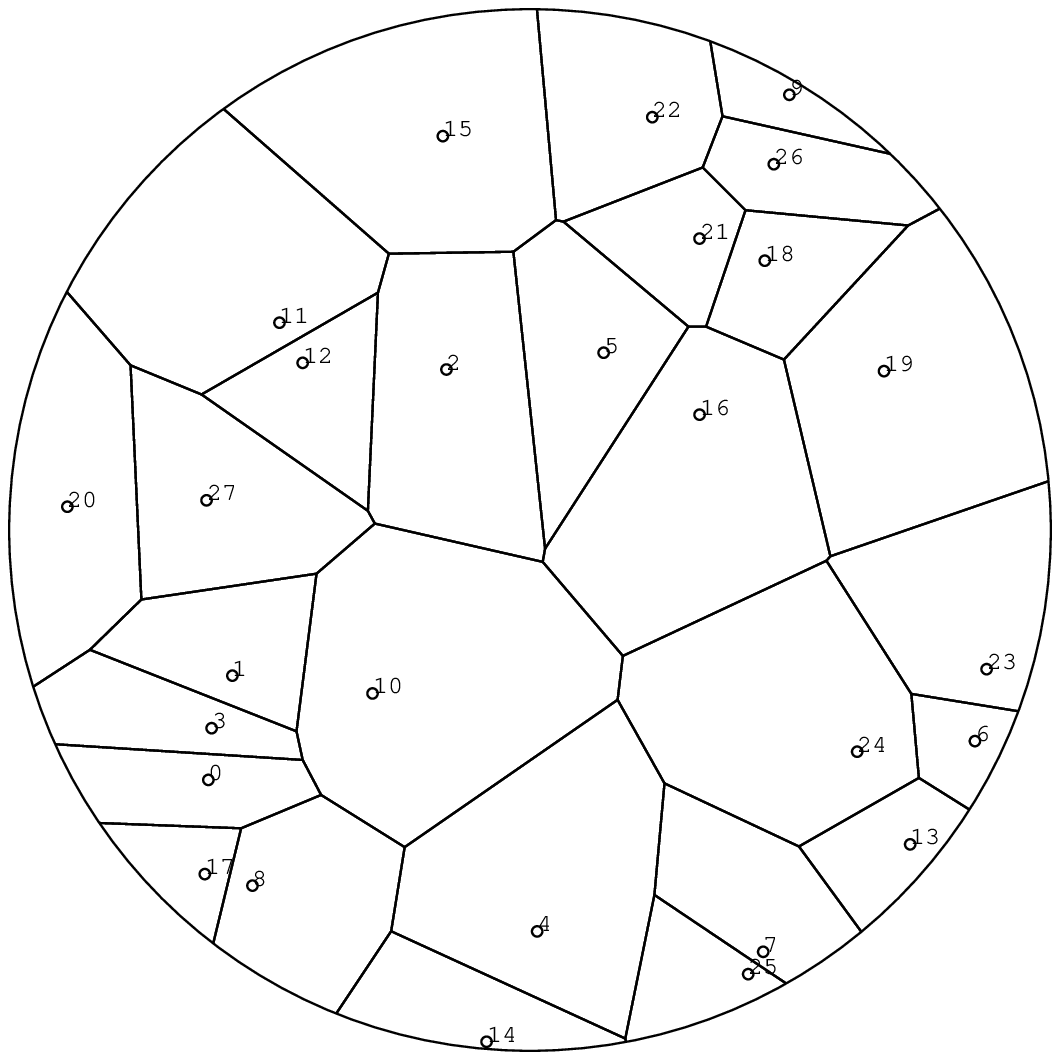}
\caption{
Illustration of an initial, randomized placement
of $28$ nodes without (or prior to) iteration.
}
\label{fig.ex3}
\end{figure}

Examples of nodes and associated tiles on $N=29$ nodes
are plotted in Figures \ref{fig.ex1} and \ref{fig.ex2}\@.
Since $29$ is a prime number, the usual factorization method would
need to align all points on a simple ring around the center
or place one at the center and distribute the others in 7 or 4 annuli.

\section{Implementation} 
\subsection{Specifications} 
The generic aspects of the algorithm are
\begin{itemize}
\item The definition
of weights $w_j$ as the area of tiles surrounding $N$ fixed
points. The sum is $\sum_j w_j=\pi$, the area of the unit circle.
\item
The specification of these tiles via iterated splitting 
of the full area with line sections that cut
mid-way between any pair of two points. This is
the Voronoi tessellation, and the definition of
Brioullin zones in solid state physics.
\item
The re-definition of points via tiles as their center-of-mass
coordinates (centroidal tessellation \cite{DuSIAM41,DuAMC133}).
\end{itemize}
The computation of polygonal areas poses no further problem, and
is done by adding the triangular areas of the sides as viewed from
their centers. The points near the rim of the circle are individually
detected and the sub-area of their spherical caps is included in
their weight.

\subsection{Modes of Use} 

The program has a mode of ``fix-point'' use
in which the number and coordinates of nodes are preassigned
and read in from an ASCII file. The program just subdivides the region
and computes the weights (areas) if no relocation loop is requested.

The other invocation uses a random-number generator
to initialize point positions in the
circle area. With a variable seed of the random-number generator,
different point sets can be initiated. 

In both cases, a number (defaulting to 0) of iterations
of moving the points to the tile centers and adjusting
the polygons may be specified. This achieves
``self-consistent'' stable results with small variances
in the weights after a modest number of loops.
If that loop count is kept at zero, the points are left frozen
at their
initial positions, and shapes like
in Figure \ref{fig.ex3} with a larger spread of tile areas result.

The generic output is a list of the two Cartesian (optionally 
two circular)
coordinates and one weight of each of the $N$ points
in ASCII format, optionally bracketed
for inclusion as arrays in other programing codes.

As a visual aid, PostScript images showing the nodes and tiles
may be created, as seen in the figures above.

\section{Summary}
We have implemented a strategy of dividing the interior
of the unit circle into distinct non-overlapping tiles
defined by halving the distance between neighboring
abscissa (nodal) points. These tiles define
weights and nodes for Riemann integration over the
unit circle.
The C++ source code is made available in the ancillary
files.

\appendix
\section{Installation}
The program is compiled by moving to the source directory
and calling the C++ compiler with

\begin{verb}
make
\end{verb}

and moving the executable \texttt{RiemCirc} to a directory
where the operating system will find it---depending on definitions
of environment variables. Alternatively the GNU autotools compile this with

\verb+ autoreconf -i -s ;+ 

\verb+ ./configure --prefix=+\textit{directory}\verb+ ;+

\verb+ make+

The Unix man-page
\texttt{RiemCirc.1} may also be moved to the subdirectory
\verb+man/man1+ of the standard folder locations or be read
with

\verb+nroff -man RiemCirc.1 | more+

In case that \texttt{doxygen} is available, one may
also call
\begin{verb}
make doc
\end{verb}
to generate a PDF file with the API\@.

\section{Synopsis}
A summary of the two versions of calling the
program is

\verb+RiemCirc+
[\verb+-v+]
[\verb+-C+]
\verb+-i+ \textit{file.asc}
[\verb+-p+ \textit{file.ps}]
[\verb+-a+]
[\verb+-l+ \textit{Nrelo}]

\verb+RiemCirc+
[\verb+-v+]
[\verb+-C+]
[\verb+-p+ \textit{file.ps}]
[\verb+-a+]
[\verb+-l+ \textit{Nrelo}]
[\verb+-r+ \textit{seed}]
[\verb+-L+ \textit{Nitr}]
\textit{N}

Brackets surround optional arguments.

The first variant reads the initial positions from \textit{file.asc},
whereas the second variant creates them randomly (if the \textit{seed}
is specified) or in annuli (if the \text{seed} is not used).

The option \verb+-v+ causes the variance of the areas to be
reported.

The option \verb+-C+ means that C-style curly parentheses will be
added to the output.

The option \verb+-p+ lets the program generate the PostScript file \textit{file.ps}
in the style seen in the manuscript.

The option \verb+-a+ causes an output of point locations
in $r$ and $\varphi$ coordinates 
rather than the default $x$ and $y$ coordinates.

The option \verb+-l+ lets the program run \textit{Nrelo}
relocation and re-tessellation loops on the initial placements before
creating the outputs.

The option \verb+-r+ places the $N$ points initially at random inside
the unit circle. The result will look more like Figure \ref{fig.ex1}
than the ring oriented Figure \ref{fig.ex2} after some iterations.

The option \verb+-L+ runs the iterations for \textit{Nitr}
different randomized starting positions and reports only
the set of positions that created the minimum variance
in the areas (weights) attributed to the relocated point sets.
So finding well-balanced point sets means using \texttt{-r}
in conjunction with suitably large \texttt{-l} and \texttt{-L}

\section{Examples}

The ancillary directory contains results of this strategy
for all $N$ points in the range $4\le N<20$ and also for some larger 
3-smooth $N$; they are generated with \texttt{make check}. 
Each result is presented in (i) a PostScript file plus (ii) an ASCII file
with a list of $(x,y)$ coordinates and their weights.
The PostScript file has a very simple format and the coordinates
of the two terminal points of each line segment can be tabulated
by selecting the \texttt{moveto} lines and scaling all 4 numbers
by 250:

\verb+grep moveto+ \textit{file}\verb+.ps | grev -v show+
\verb+| awk '{printf($1/250,$2/250,$4/250,$5/250)}'+

\bibliography{all}

\end{document}